\documentclass[letterpaper, 10 pt, conference]{ieeeconf}  

\IEEEoverridecommandlockouts                              

\usepackage{cite}
\usepackage{amsmath,amssymb,amsfonts}
\usepackage{algorithmic}
\usepackage{graphicx}

\usepackage{subfigure}

\usepackage{textcomp}
\usepackage{xcolor}
\usepackage{comment}
\def\BibTeX{{\rm B\kern-.05em{\sc i\kern-.025em b}\kern-.08em
    T\kern-.1667em\lower.7ex\hbox{E}\kern-.125emX}}

\usepackage[utf8]{inputenc}
\usepackage{amsmath}
\usepackage{amsfonts}
\usepackage{amssymb}

\usepackage{times}
\usepackage{mathtools}
\usepackage{calrsfs}
\usepackage{comment}
\usepackage{float}

\DeclareMathAlphabet{\pazocal}{OMS}{zplm}{m}{n}

\usepackage[ruled,longend]{algorithm2e}
\newcommand{\tr}[1]{{#1}^\mathsf{T}}

\usepackage{algorithmic}
\usepackage{array}
\usepackage{bm}    

\usepackage{braket} 

\usepackage[export]{adjustbox}

\title{\LARGE \bf  Receding Horizon Differential Dynamic Programming Under Parametric Uncertainty
}
\author{Yuichiro Aoyama$^{1,2}$, Augustinos D. Saravanos$^{1}$ and Evangelos A. Theodorou$^{1}$
\thanks{$^{1}$School of Aerospace Engineering, Georgia Institute of Technology, Atlanta, GA, USA}
\thanks{$^{2}$Komatsu Ltd.,Tokyo, Japan}
\thanks{\{yaoyama3, asaravanos3, evangelos.theodorou\}@gatech.edu}}

\begin{document}
\maketitle
\begin{abstract}
Generalized Polynomial Chaos (gPC) theory has been widely used for representing parametric uncertainty in a system, thanks to its ability to propagate uncertainty evolution. In an optimal control context, gPC can be combined with several optimization techniques to achieve a control policy that handles effectively this type of uncertainty. Such a suitable method is Differential Dynamic Programming (DDP), leading to an algorithm that inherits the scalability to high-dimensional systems and fast convergence nature of the latter. In this paper, we expand this combination aiming to acquire probabilistic guarantees on the satisfaction of nonlinear constraints. In particular, we exploit the ability of gPC to express higher order moments of the uncertainty distribution - without any Gaussianity assumption - and we incorporate chance constraints that lead to expressions involving the state covariance. Furthermore, we demonstrate that by implementing our algorithm in a receding horizon fashion, we are able to compute control policies that effectively reduce the accumulation of uncertainty on the trajectory. The applicability of our method is verified through simulation results on a differential wheeled robot and a quadrotor that perform obstacle avoidance tasks.

\end{abstract}


\section{Introduction}




One of the most challenging problems in the control field arises when a system operates under uncertainty.
To address this issue, existing approaches can be classified into model-free and model-based ones.
The former class has found several successful applications in the context of reinforcement learning (e.g., \cite{kober2013reinforcement, theodorou2010generalized, buchli2011learning}), however main drawbacks of such methods include their requirement for numerous interactions with the physical system and slow convergence rates. 
On the other hand, methods that belong in the latter category (e.g., \cite{ddp, li2004iterative, lantoine2012hybrid}) can be significantly faster, but their performance relies substantially on the accuracy of the model. Moreover, most of them, such as \cite{todorov2005generalized} and \cite{Theodorou2010SDDP}, utilize stochastic differential equations for representing uncertainty with Brownian motion which assumes Gaussianity. In addition, these methods are not taking into account any stochasticity arising from internal model parameters.

To address this problem, generalized Polynomial Chaos (gPC) theory \cite{Xiu2002askey_pc} has been used for representing parametric uncertainty. This approach approximates stochastic processes using orthogonal polynomials which are chosen according to the type of distribution of the uncertain parameters. Therefore, a remarkable advantage of gPC is that it does not require any Gaussian assumptions on the dynamics. Moreover, by combining gPC with dynamical systems, it becomes possible to propagate uncertainty evolution and express higher order moments of its distribution.

From an optimal control perspective, gPC has found some applications in controlling systems under parametric uncertainty. In \cite{Fagiano2012nonlin_mpc_pce}, a Model Predictive Control (MPC) with gPC approach is proposed, taking into account the expectations of constraints. Another gPC-MPC method is suggested in \cite{Kim2013gpc_approx_mpc} for dealing with additive Gaussian noise, while satisfying linear chance-constraints. In \cite{Andrian2019}, the authors developed an offset-free MPC based on gPC under linear constraints. By linearizing the constraints, \cite{Lucia2015efficient_mpc_pce} presents an effective way to address constrained MPC with gPC.
Furthermore, an efficient method of solving the gPC-MPC problem was demonstrated in \cite{Mesbah2014stoch_MPC_prob_const}, by exploiting the second-order cone constraints that emerge from chance constraints. Finally, a chance-constrained stochastic nonlinear control approach with gPC was presented in \cite{Nakka2019chance_constraints}, demonstrating a robot performing a reaching task while dodging an spherical obstacle.

Differential Dynamic Programming (DDP) is a trajectory optimization method that was first introduced in \cite{ddp} and has found several successful applications such as \cite{kumar2016optimal, morimoto2003minimax, boutselis_ddp_lie, ypan_pddp_journal} etc. 
Its main advantages include its quadratic convergence rate and its greater scalability to high-dimensional systems compared to other optimization techniques. 
For this reason, combining DDP with gPC \cite{George2019stochastic_mechsystem} can lead to an algorithm that maintains these attributes, while being able to handle parametric uncertainty effectively.

In this paper, we extend the algorithmic framework suggested in \cite{George2019stochastic_mechsystem}, by proposing a receding horizon method that handles effectively parametric uncertainty and nonlinear constraints while maintaining proper convergence properties and scalability to high-dimensional systems. 
In particular, by taking advantage of our recent results on constrained DDP \cite{aoyama2020constrained}, we incorporate nonlinear chance constraints that through gPC theory yield deterministic constraints that involve the state covariance of the system.
Moreover, in order to deal with the growth of uncertainty through the trajectory, we integrate feedback by executing our algorithm in a receding horizon scheme. As confirmed by the simulation results, our method can satisfy multiple nonlinear constraints under parametric uncertainty, while also being applicable to high-dimensional robotic tasks.

The remaining of the paper is organized as follows. A brief overview of gPC theory and DDP is provided in Section II. In Section III, we propose a receding horizon approach that combines DDP and gPC while satisfying nonlinear chance constraints. In Section IV, we verify the effectiveness of the suggested algorithm through simulation results on a differential wheeled robot and a quadrotor. The conclusions of our work along with future directions are provided in Section V.

\section{Preliminaries}
In this section, we present some necessary preliminaries. We start with covering the theory of gPC and showing how it can represent dynamical systems with parametric uncertainty. Next, we give a concise overview of DDP and we describe how it can be applied to dynamical systems represented by gPC.  
\label{sec:preliminaries}
\subsection{Generalized Polynomial Chaos}
\label{sec:gpc}
\subsubsection{Polynomial Chaos Expansion}
Let us define a probability space $(\Omega,\pazocal{F},P)$ where $\Omega$ is the sample space, $\pazocal{F}$ is the $\sigma$-field of $\Omega$, and $P$ is the probability measure. Let $\omega \in \Omega$ and ${\xi}(\omega) = ( \xi_{1}(\omega), \dots, \xi_{d}(\omega) ) \in \mathbb{R}^d$ be a continuous random variable vector with mutually independent components whose distributions are subject to the probability density function (PDF) $\rho(\xi)$ with support $I_{\xi}$. Furthermore, let $F(\xi)$ be a function of the random variable $\xi$ and $L^{2}_{\rho}$ be the space of square integrable functions with norm $||F||_{L^2_{\rho}} = \mathbb{E}[F^2]^{1/2}$, i.e.:
\[ L^{2}_{\rho} = \Set{F:I_{\xi} \rightarrow \mathbb{R} |
\begin{array}{l}
\int_{I_{\xi}}F^{2}(\xi)\rho(\xi)d\xi \rightarrow  \infty\
\end{array}}. \]
Then, the polynomial chaos expansion of $F$ can be obtained as \cite{Xiu2010gpc_text}:
\begin{align}
F(\xi)\approx \sum_{j = 0}^{\infty} F_{j}\Phi_{j}(\xi). \label{gPC expansion}
\end{align}
The functions $\Phi_{j}(\xi)$ in expression (\ref{gPC expansion}) are the generalized polynomial chaos basis functions, that satisfy orthogonality:
\begin{align}\label{eq:exp_phim_phi_n}
    \mathbb{E}[\Phi_m(\xi)\Phi_n(\xi)] = \int_{I_{\xi}}\Phi_m(\xi)\Phi_n(\xi)\rho(\xi)d\xi = \gamma \delta_{mn}, \nonumber
\end{align}
where $\gamma = \mathbb{E}[\Phi^2(\xi)] = \langle \Phi_{j} \Phi_{j}\rangle$ are the normalization factors and $\delta_{mn} = 0$ if $m \neq n$, $\delta_{mn} = 1$ if $m = n$.  

For several types of distributions, specific gPC basis functions exist for describing the underlying random variables of a given function. Some of them are listed in Table \ref{tab:polynomials}.


In practical applications, a truncated version of the gPC expansion (\ref{gPC expansion}) can be used by taking into account polynomials of degree up to $r$:
\begin{align*}
    F(\xi)\approx \sum_{j = 0}^{K} F_{j}\Phi_{j}(\xi),
\quad
K = \frac{(r+d)!}{r!d!} - 1,
\end{align*}
where $K$ is the total number of polynomials.
The coefficients of this expansion can be obtained by exploiting the orthogonality of ${\Phi}$. Taking the inner product of $F$ and $\Phi_{j}$ yields:
\begin{align*}
    F_{j} = \frac{\int_{I_{\xi}}F(\xi)\Phi_{j}(\xi)\rho(\xi)d{\xi}}{\int_{I_{\xi}}\Phi^{2}_{j}(\xi)\rho(\xi)d{\xi}},
\end{align*}
which is called Galerkin projection.
Note that the gPC basis functions $\Phi_{j}(\xi)$ are the products of the polynomials of mutually independent random variables with multi-index $|i| = i_1 + \dots + i_d$. More precisely, we have:
\begin{align*}
\Phi_{i}(\xi) = \phi_{i1}(\xi_{1})\dots\phi_{id}(\xi_{d}),\quad 0\leq |i| \leq r.
\end{align*}

\subsubsection{Dynamical systems with polynomial chaos}
In order to apply DDP to a stochastic system represented by gPC, information about the evolution of the system including the dynamics and the Jacobian matrix is necessary. 
We consider the dynamical system whose state $x(t) \in\mathbb{R}^{n}$ and initial state $x(t_{0})=x_{0}$ are influenced by the uncertain parameters $\zeta^{p} \in{\mathbb{R}}^{d_{p}}$ and $\zeta^{0} \in{\mathbb{R}}^{d_{0}}$, respectively. It is also assumed that the distributions of these parameters are known. The control $u(t) \in\mathbb{R}^{m}$ is set to be deterministic. The system dynamics 
$f:\mathbb{R}^{n}\times \mathbb{R}^{m}\rightarrow\mathbb{R}^{n}$ are provided by:
\begin{align}
\dot{x}(t) = f(x(t), u(t), t;\zeta^{p}),\quad
x(t_{0}) = x(t_{0};\zeta^{0}).
\end{align}
Applying the gPC expansion to random variables $\zeta^{p}$ and $\zeta^{0}$ with respect to the mutually independent standard random variables $\xi^{p}$ and $\xi^{0}$ yields:
\begin{align*}
\zeta_{i}^{p} \approx \sum_{j = 0}^{Z}\zeta_{ij}^{p}\phi_{j}^{p}(\xi^{p}), \quad
\zeta_{i}^{0} \approx \sum_{j = 0}^{Z_0}\zeta_{ij}^{0}\phi_{j}^{0}(\xi^{0}),
\end{align*}
where $Z$ and $Z_{0}$ are the orders of the polynomials.
Likewise, the initial states are expanded as:
\begin{align*}
x_i(t_0) \approx \sum_{j=0}^{K_{0}}x_{ij}(t_{0})\Phi_{j}^{0}(\xi^{0}).
\end{align*}

\begin{table}
\centering
\caption{Correspondence between the type of underlying random variables and the type of gpc basis functions.}
\begin{tabular}{ ccc } 
 Distribution & PDF: $\rho$ & gPC basis polynomials \\
  \hline
Gaussian &  $\frac{1}{\sqrt{2\pi}}e^{-x^{2}/2}$ & Hermite\\
Beta &$\frac{(1-x)^{\alpha}(1-x)^{\beta}}{2^{\alpha + \beta +1} B (\alpha+1, \beta +1)}$ & Jacobi\\
Uniform &$\frac{1}{2}$ & Legendre\\
Gamma &$\frac{x^{\alpha}{e^{-x}}}{\Gamma(\alpha + 1)}$ & Laguerre\\
 \hline
\end{tabular}
\label{tab:polynomials}
\end{table}

In order to parameterize $x$ with both $\xi^p$ and $\xi^0$, we define the concatenated vector of random parameters as $\xi = ( \xi^p,\xi^0 ) \in  \mathbb{R}^d$, with $d = d_{p} +d_{0}$. 
Since the evolution of $x$ depends on $\xi$, the gPC expansion of the state vector and its time derivatives are given by:
\begin{align*}
{x}_{i}(t) \approx \sum_{j=0}^{K}{x}_{ij}(t)\Phi_{j}(\xi), \quad
\dot{x}_{i}(t) \approx \sum_{j=0}^{K}\dot{x}_{ij}(t)\Phi_{j}(\xi).
\end{align*}
Applying Galerkin projection to these equations yields the coefficients:
\begin{align*}
    x_{ij} = \frac{\int_{I_{\xi}}x_{i}(t)\Phi_{j}(\xi)d{\xi}}{{\langle  \Phi_{j}, \Phi_{j}\rangle }}, \quad
    \dot{x}_{ij} = \frac{\int_{I_{\xi}}\dot{x}_{i}(t)\Phi_{j}(\xi)d{\xi}}{{\langle  \Phi_{j}, \Phi_{j}\rangle }}.
\end{align*}
%
Therefore, the evolution of the system can be obtained as an evolution of gPC coefficients, i.e.:
\begin{align*}
    \dot{X}(t) &= f(X(t), u(t),t),\\
    \text{where} \quad X &= (x_{10},\dots,x_{1K},\dots,x_{nK})\in{\mathbb{R}}^{n(K+1)}.
\end{align*}
%
%
Finally, differentiating the time derivative of the coefficients with respect to the coefficients of $x$ and $u$ yields the elements of the Jacobian matrix:
\begin{align*}
\frac{{\partial}\dot{x}_{ij}}{\partial{x}_{gh}}
&= \frac{\int_{I_{\xi}} \frac{{\partial}{f_{i}}}{\partial{x_{gh}}}\Phi_{j}(\xi)\rho(\xi)d\xi}{\langle  \Phi_{j}, \Phi_{j}\rangle }
= \frac{\int_{I_{\xi}} \frac{{\partial}{f_{i}}}{\partial{x_{g}}} \Phi_{h}(\xi) \Phi_{j}(\xi)\rho(\xi)d\xi}{\langle  \Phi_{j}, \Phi_{j}\rangle },
\\
\frac{{\partial}\dot{x}_{ij}}{\partial{u}_{k}}
& = \frac{\int_{I_{\xi}} \frac{{\partial}{f_{i}}}{\partial{u}}\Phi_{j}(\xi)\rho(\xi)d\xi}{\langle  \Phi_{j}, \Phi_{j}\rangle },\\  
\text {with}& \quad
i,g = 1,\dots, n,\quad j,h = 0,\dots, K,\quad k = 1,\dots, m.
\end{align*}
The integrals can be numerically evaluated online by Gaussian quadrature. The nodes and weights for the computation can be obtained by solving an eigenvalue problem \cite{Golb1969gaussquad}.

Once the gPC expansion of $x_{i}(\xi)$ is obtained, statistical information can be computed analytically. The mean, variance, and covariance are given by:
\begin{align*}
     \mathbb{E}[x_{i}(t)] &= \int_{I_{\xi}}\sum_{i=0}^{K}x_{ij}(t)\Phi_{j}(\xi)\rho d\xi = x_{i0}(t),\\
     \mathbb{VAR}[x_{i}(t)] 
     &= \sum_{j=1}^{K}x_{ij}^{2}(t)\langle  \Phi_{j},\Phi_{j}\rangle,\\
     \mathbb{COV}[x_{i}(t), x_{g}(t)] &= \sum_{j=1}^{K}x_{ij}(t)x_{gj}(t)\langle  \Phi_{j}, \Phi_{j}\rangle.
 \end{align*}
%
%

\subsection{Differential Dynamic Programming}
\label{subsec:ddp}
\subsubsection{Deterministic DDP}
Next, we provide a brief overview of discrete-time DDP. For more details the reader can be referred to \cite{ddp, li2004iterative}.
Let us consider the discrete-time finite-horizon optimal control problem:
\begin{align}
\label{unconstrained_optimalcontrol}
& \min_{\bm{U}}\hspace{0.8mm}J(\bm{X}_{\rm{d}},\bm{U})=\min_{\bm{U}}\big[{\textstyle\sum_{k=0}^{N-1} l(x_k,u_k)}+\phi(x_N)\big]\\\notag
& \ \text{subject to} \quad x_{k+1}=f(x_k,u_k),\, \ k=0, \dots ,N-1,
\end{align}
where $x_{k} = x(t_k)$, $u_k = u(t_k)$
and $f$ is the transition dynamics function. The scalar functions $l$, $\phi$ and $J$ correspond to the running, terminal and total cost, respectively. With ${\bm{X}}_{\rm{d}} :=({\tr{x}_{0}},\dots,{\tr{x}_{N}})$ and ${\bm{U}} :=({\tr{u}_{0}},\dots,{\tr{u}_{N-1}})$, we denote the state and control sequences, respectively, over the time horizon.

Next, the notion of the value function is introduced: 
\begin{align*}\label{value-function-definition}
V_k(x_k):=\min_{u_k}J(\bm{X}_{\rm{d}},\bm{U}),
\end{align*}
which provides the minimum cost-to-go at each state and time.
The Bellman's principle of optimality can be expressed through the following backward propagation rule:
\begin{equation}\label{bellman}
    V_k(x_k) =  \min_{u_k} [l(x_k,u_k) + V_{k+1}(x_{k+1})].
\end{equation}

Let us also define the following function: 
\begin{equation*}
\label{Qfunction}
Q_k(x_k, u_k) := l(x_k, u_k) + V_{k+1}(x_{k+1}).
\end{equation*}
which is the quantity to be minimized in \eqref{bellman}.
During the backward pass of DDP, problem \eqref{unconstrained_optimalcontrol} is solved locally by expanding both sides of \eqref{bellman} about some given nominal trajectories $\bar{\bm{ X}}_{\rm{d}}$, $\bar{\bm{ U}}$. 
By taking the quadratic expansions of $Q_k$ around the vicinity of $\bar{\bm{ X}}_{\rm{d}}$, $\bar{\bm{ U}}$, we get:
\begin{align*}
  &Q_{k}(x_k,u_k)\approx  Q_{k}+{Q}_{x,k}^{\mathsf{T}}\delta x_{k}+{Q}_{u,k}^{\mathsf{T}}\delta u_{k}+\\
  &\quad\quad\textstyle{\frac{1}{2}}(\tr{\delta x}_k{Q}_{xx,k}\delta x_{k}+2\tr{\delta x}_k{Q}_{xu,k}\delta u_{k}+\tr{\delta u}_k{Q}_{uu,k}{\delta}{u}_k),
\end{align*}
where:
%
%
\begin{equation}
\label{Qexpanded}
\begin{split}
    &{Q}_{xx,k} ={l}_{xx}+f_x^{\mathsf{T}}{V}_{xx,k+1} f_x,\, \ {Q}_{x,k}={l}_{x}+f_x^{\mathsf{T}}{V}_{x,k+1}
     \\
    &{Q}_{uu,k} = {l_{uu}}+f_u^{\mathsf{T}}{{V}_{xx,k+1}} f_u, \, \ {Q}_{u,k}={l_{ u}}+f_u^{\mathsf{T}}{{V}_{x,k+1}}
    \\
    &{Q}_{xu,k} ={l_{xu}}+f_x^{\mathsf{T}}{{V}_{xx,k+1}} f_u.
\end{split}
\end{equation}
with $\delta x_k:= x_k -\bar{x}_k$, $\delta u_k:=u_k - \bar{u}_k$ being the deviations about the nominal sequences. Note that in \eqref{Qexpanded}, the $Q$ functions are evaluated on $\bar{\bm{ X}}_{\rm{d}}$, $\bar{\bm{ U}}$. Given the expressions \eqref{Qexpanded}, we can now explicitly minimize \eqref{bellman} with respect to $\delta\bm{u}$ in order to acquire the locally optimal control deviations as:
\begin{equation*}
\label{delta-u-star}
\begin{split}
    & \delta u^{\ast}_k=\bm{k}_k+\bm{K}_k\delta x_{k},\\
    \text{with}\quad &\bm{k} := -{Q}^{-1}_{uu}{Q_{u}},\hspace{1.8mm} \bm{K} := -{Q}^{-1}_{uu}{Q_{ux}}.
\end{split}
\end{equation*}
Moreover, the value function $V_k$ is also quadratically expanded and by plugging it into \eqref{bellman} along with $\delta u^\ast_k$, we get:
\begin{equation} \label{eq:value_riccati}
\begin{split}
    V_{x,k}&={Q_{x,k}}-{Q_{xu,k}}{Q_{uu,k}^{-1}}{Q_{u,k}}\\
    V_{xx,k}&={Q_{xx,k}}-{Q_{xu,k}}Q_{uu,k}^{-1}{Q_{ux,k}}.
    \end{split}
\end{equation}
The equations \eqref{eq:value_riccati} are propagated backwards in time using the terminal condition $V(x_N) = \phi(x_{N})$. 

Subsequently, in the forward pass, the new control sequence is applied to the system.
The resulting trajectories will be used as the nominal ones at the next backward pass leading to an iterative process, which is terminated with the satisfaction of some predefined convergence criteria.
Regarding the cost functions, control effort and deviation from the desired state $x^{\rm{d}}$ are typically penalized, i.e.:
\begin{align*}
     l(x_{k},u_{k}) &= \frac{1}{2}(x_{k}-x^{\rm{d}}_k)^{\mathsf{T}}A(x_{k}-x^{\rm{d}}_{k})+\frac{1}{2}u_{k}^{\mathsf{T}}Ru_{k}, \label{eq:running_cost}\\
    \phi(x_{N}) &= \frac{1}{2}(x_{N}-x_{N}^{\rm{d}})^{\mathsf{T}}A_{\rm{f}}(x_{N}-x_{N}^{\rm{d}}),
\end{align*}
where $A_{\rm{f}}$, $R$ are positive definite matrices and $A$ is a positive semidefinite matrix.
\subsubsection{DDP with gPC}
We will now demonstrate how unconstrained DDP can be combined with gPC. Since the dynamics are required to be discretized for applying discrete-time DDP, we can apply Euler discretization on them.
Moreover, in order to obtain the optimal control of a system represented by gPC, the expectation of the cost can be used \cite{Fisher2011optimal_trj_gpc}.
 Taking the expectation of the running cost of the discretized system yields:
\begin{align*}
    \mathbb{E}[l(X_{k},u_{k})] &= \frac{1}{2}(X_{k}-X^{\rm{d}}_{k})^{\mathsf{T}}A^{X}(X_{k}-X^{\rm{d}}_{k}) + \frac{1}{2}u_{k}^{\mathsf{T}}Ru_{k},\\
\text{with}\quad
X^{\rm{d}}_{k} &= (x^{\rm{d}}_{k,1},0_{1\times K},
x^{\rm{d}}_{k,2},0_{1\times K},\dots,\\
& \quad \quad \quad x^{\rm{d}}_{k,n},0_{1\times K})^{\mathsf{T}}\in{\mathbb{R}}^{n(K+1)},\\
A^{X} &= A \otimes {\rm{diag}}(\langle \Phi_{0},\Phi_{0} \rangle, \langle \Phi_{1},\Phi_{1} \rangle, \dots \langle \Phi_{K},\Phi_{K} \rangle),
\end{align*}
where $X_k^{d}$ is the expanded desired state $x_{d}$ and $\otimes$ stands for the Kronecker product. Note that since $x_{d}$ is deterministic, $X_{d}$ only contains terms corresponding to the mean of $x_{d}$ and zeros for higher order moments.
The expectation of the terminal cost can be simplified using the same procedure.
In order to penalize not only the mean but also the higher order moments of the state, $A$ can be modified by taking into account elements corresponding to the higher order moments ($x_{ij}, i \neq 0$) as: 
\begin{align*}
    A^{X} = \begin{bmatrix}
    A_{1} &\dots& 0_{K+1}\\
    \vdots& \ddots &\vdots\\
    0_{K+1}& \cdots& A_{n}
    \end{bmatrix}
\end{align*}
with 
$A_{i} = {\rm{diag}}(a_{i0}, a_{i1}\langle \Phi_{1},\Phi_{1} \rangle,\dots,a_{iK}\langle \Phi_{K},\Phi_{K} \rangle )$. More details on the cost derivation can be found in \cite{George2019stochastic_mechsystem}.

\section{Receding Horizon Constrained DDP under Parametric Uncertainty} \label{sec:const gpc DDP}
\subsection{Problem Formulation}
Let us consider the following finite-horizon optimal control problem under parametric uncertainty and subject to control and state constraints:
\begin{subequations} \label{chanceconstraint_optimalcontrol}
\begin{align}
&\min_{\bm{U}}\hspace{0.8mm}J(\bm{X}_{\rm{d}},\bm{U})=\min_{\bm{U}}\big[{\textstyle\sum_{k=0}^{N-1} l(x_k,u_k)}+\phi(x_N)\big]\\
&\text{subject to}\quad x_{k+1}=f(x_k,u_k;\zeta), \\
&\qquad \qquad \quad {\rm{P}}[g(x_{k},u_{k})\leq 0] > p_{\rm{c}},\label{eq:chanceconstraint_optimalcontrol_const}\\
&\qquad \qquad \quad u_{\rm{min}} \leq u_{k} \leq u_{\rm{max}}, \label{eq:chanceconstraint_optimalcontrol_contollim} \\
&\qquad \qquad \quad k=0,\dots,N-1. \notag
\end{align}
\end{subequations}
%
%
The chance constraint \eqref{eq:chanceconstraint_optimalcontrol_const} indicates that $g(x_{k},u_{k})\leq 0$ is satisfied with a probability greater than 
a specified value $p_{c}$. 
By taking the gPC expansion of the dynamics and the expectation of the performance index, problem \eqref{chanceconstraint_optimalcontrol} can be transformed into the following deterministic one:
\begin{subequations} \label{chanceconstraint_optimalcontrol_gpc} 
\begin{align}
& \min_{\bm{U}} \hspace{0.8mm}J(\bm{X},\bm{U})=\min_{\bm{U}}\mathbb{E}\big[{\textstyle\sum_{k=0}^{N-1} l(X_k,u_k)}+\phi(X_N)\big]\\
&\text{subject to}\quad X_{k+1}=f(X_k,u_k), \\
&\qquad \qquad \quad g_{\rm{g}}(X_{k},u_{k}) < \alpha,  \label{eq:deterministicconstraint_optimalconrtol_gpc_const} \\
&\qquad \qquad \quad u_{\rm{min}} \leq u_{k} \leq u_{\rm{max}}, \\
&\qquad \qquad \quad k=0, \dots ,N-1, \notag
\end{align}
\end{subequations}
%
%
%
where ${\bm{X}} :=({\tr{X}_{0}},\dots,{\tr{X}_{N}})$. The chance constraint \eqref{eq:chanceconstraint_optimalcontrol_const} has now been transformed into the deterministic \eqref{eq:deterministicconstraint_optimalconrtol_gpc_const} by using the appropriate function $g_{\rm{g}}$ and constant $\alpha$. The choice of these two is further explained during the simulation examples in Section \ref{sec:results}.

\subsection{Constrained DDP with generalized Polynomial Chaos}
After expressing the originally stochastic problem \eqref{chanceconstraint_optimalcontrol} in a deterministic manner, we are now able to 
solve it using a deterministic constrained trajectory optimization method such as Augmented Lagrangian (AL) DDP \cite{aoyama2020constrained}. 
From now on, we will refer to the resulting approach as gPC Constrained DDP (gPC CDDP). 
The constraints \eqref{eq:deterministicconstraint_optimalconrtol_gpc_const} are first, redefined as $G = g_{\rm{g}}-\alpha < 0$, where $G = (G_{1},\dots,G_{w})^{\mathsf{T}}$ is a vector of $w$ constraints, and then, integrated into the original objective function by adding a function $\pazocal {P}$ which penalizes their violation. Therefore, problem \eqref{chanceconstraint_optimalcontrol_gpc} can now be expressed as:
\begin{subequations} \label{eq:constrained_control_multipliers}
\begin{align}
& \min_{\bm U}\hspace{1mm} \Big[
J( \bm X, \bm U) + \sum_{i=1}^{w}\sum_{k=0}^{N-1} \pazocal{P}(\lambda_i^k, \mu_i^k, G_{i,k}(X_k,u_k))
\Big] \\
& \text{subject to}\quad X_{k+1}=f(X_k, u_k),  \\
& \qquad \qquad \quad u_{\rm{min}} \leq u_{k} \leq u_{\rm{max}}, \\
& \qquad \qquad \quad k = 0, \dots, N-1, \notag \label{eq:constrained_control_multiplier_contollim}
\end{align}
\end{subequations}
where $\lambda$ and $\mu$ are the Lagrange multipliers and penalty parameters, respectively. It is known that \eqref{eq:constrained_control_multipliers} can provide a solution to \eqref{chanceconstraint_optimalcontrol_gpc} under mild assumptions \cite{birgin2005numerical}.

The optimization process of AL DDP consists of two stages: an inner and an outer loop. In the inner loop, a locally optimal solution of \eqref{eq:constrained_control_multipliers} is obtained by applying DDP. 
In the outer loop, the parameters $\lambda$ and  $\mu$ are updated according to the value of the derivative of $\pazocal{P}$ and $\pazocal{P}$ itself, respectively. As for $\mu$, they are increased monotonically, given the constraint improvement from the inner loop is not enough. More details regarding the choice of $\pazocal{P}$ can be found in \cite{birgin2005numerical}, \cite{nocedal2006numerical}.

AL DDP is quite robust as a trajectory optimization method, reaching the desired state at an early stage of the optimization iterations. However, since the constraints \eqref{eq:deterministicconstraint_optimalconrtol_gpc_const} are treated as soft constraints - as they are part of the cost - they may be slightly violated. On the other hand, control constraints can be strictly satisfied by using control-limited DDP \cite{tassa2014control} in the inner loop instead of unconstrained DDP. Thus, we use a combination of AL and control-limited DDP similar to the one presented in \cite{aoyama2020constrained}. The slight constraint violation that we mentioned can be observed especially when the horizon of the trajectory is quite long or the number of the constraints is large. In our receding horizon approach presented in \ref{sec:Receding_Horizon_DDP}, this drawback is alleviated by taking into account a shorter prediction horizon which leads to a smaller number of constraints. 

\subsection{Receding Horizon Constrained DDP} \label{sec:Receding_Horizon_DDP}
Applying gPC CDDP directly can be effective, but it may suffer from the accumulation of state variance on the trajectory. This variance growth can result in violating the chance constraints and being unable to reach the desired state target, as it is illustrated in Section \ref{sec:results}. In order to deal with this problem, we implement our proposed method in a receding horizon fashion. An similar MPC variant of DDP has been presented in \cite{tassa2007mpcddp}, but in a fully deterministic setting. In the parametric uncertainty case, the benefits of following an MPC approach are mainly the following. First of all, by taking into account the current state during the computation of each MPC control sequence, we integrate feedback into our method. Second, by considering the state evolution only during the prediction horizon $H$ - which is shorter than the whole time horizon $N$ - the growth of the variance in our computations can be maintained to lower levels. Finally, by addressing a shorter trajectory optimization problem than the original one, the constraints taken into account by gPC CDDP will be reduced.


Essentially, instead of solving \eqref{eq:constrained_control_multipliers}, we solve at each time step $k$, the following problem:
\begin{subequations} \label{eq:constrained_control_multipliers_MPC}
\begin{align}
& \min_{\bm U}\hspace{1mm} \Big[
J( \bm X, \bm U) + \sum_{i=1}^{w}\sum_{\kappa=k}^{k+H-1} \pazocal{P}(\lambda_i^{\kappa}, \mu_i^{\kappa}, G_{i,{\kappa}}(X_{\kappa},u_{\kappa}))
\Big] \\
& \text{subject to}\quad X_{{\kappa}+1}=f(X_{\kappa}, u_{\kappa}),  \\
& \qquad \qquad \quad u_{\rm{min}} \leq u_{\kappa} \leq u_{\rm{max}}, \\
& \qquad \qquad \quad \kappa = k, \dots, k+H-1. \notag
\end{align}
\end{subequations}
Our algorithm starts with computing the locally optimal control sequence for problem \eqref{eq:constrained_control_multipliers_MPC} using gPC CDDP.
Next, the first segment of the control is applied to the real system, keeping it in the confidence region of gPC CDDP while satisfying the chance constraints. Note that the specific values of the parameters are not required for the system to remain inside of the bounds. After obtaining the updated state of the real system, the gPC coefficients are updated accordingly. We can update the mean part (0-th coefficient) of the gPC dynamics by using the real system information and set the remaining coefficients - which correspond to the variance - to be zero. If we assume noisy measurements, the latter coefficients can be given non-zero values. 
Finally, in order to exploit a warm start for the next MPC computation, the control is also applied to the gPC dynamics, whose variance is low due to the information from the real system. 
This gPC MPC CDDP method can successfully keep the variance low, addressing the problem stated earlier in this section. This is further explained and demonstrated in Section \ref{sec:results}.
\section{Results} \label{sec:results}
In this section, we present simulation results that verify the effectiveness of our approach. Initially, we apply the method on a differential wheeled robot while explaining in detail its performance compared to other related approaches. Subsequently, we demonstrate the applicability of our method to more complex systems such as a quadrotor.
\subsection{Differential Wheeled Robot}\label{sec:result_robot}
\begin{figure}[t]
\centering
\includegraphics[scale=0.32]{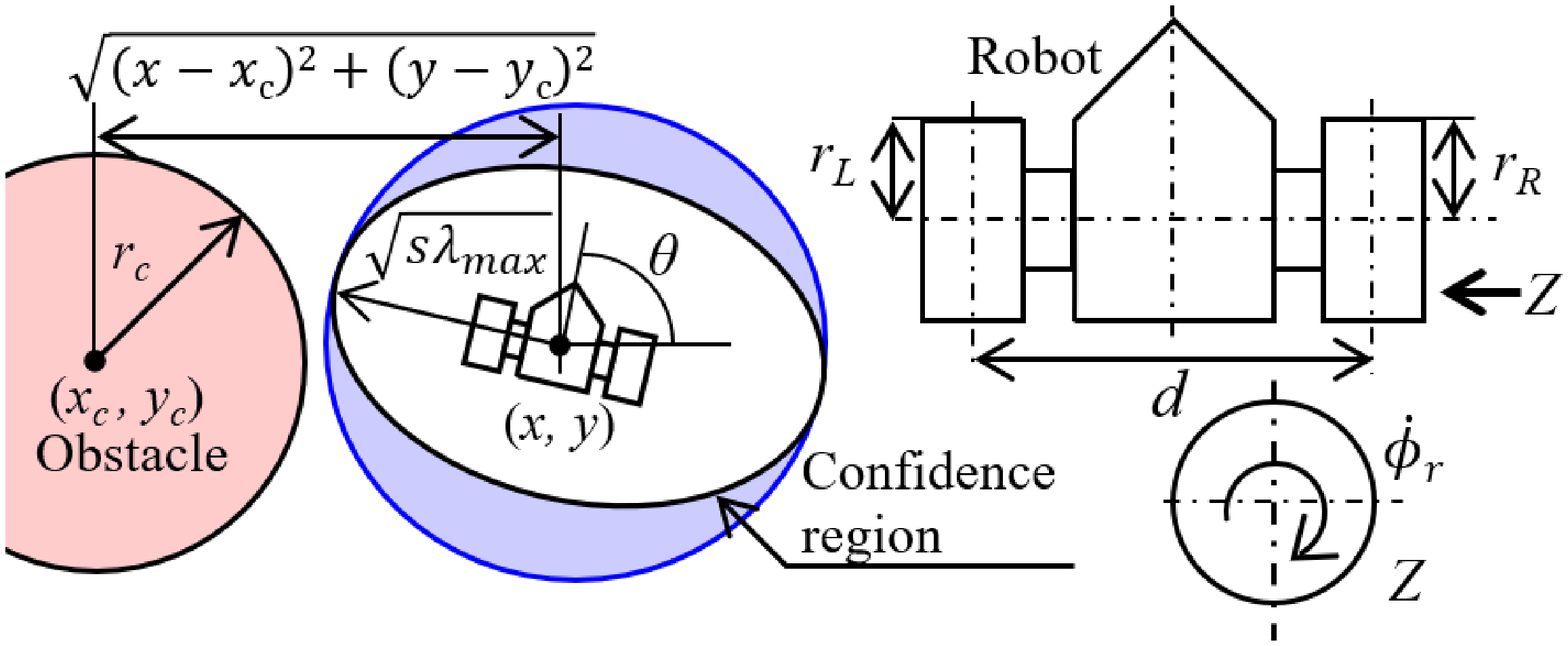}
\caption{Differential wheeled robot: an obstacle and the confidence region of the robot.}
\label{fig:robot}
\end{figure}
We will first test our algorithm on a differential wheeled robot as shown in Fig. \ref{fig:robot}. The goal of the robot is to reach a desired target while avoiding obstacles. Its model contains three uncertain parameters: the tread $d$ and the wheel radii $r_{R}$ and $r_{L}$, all of which are assumed to be normally distributed with means $0.2$ and variances $1.5\times 10^{-3}$. Hermite polynomials are chosen - according to Table \ref{tab:polynomials} - to describe these parameters. The first few of them, their inner products, and the PDF are given by:
\begin{align*}
   &\phi_{0}(z) = 1,\quad \phi_{1}(z) = z, \quad \phi_{2}(z) = z^{2}-1,\dotsc, \\ 
    &\langle \phi_{i}, \phi_{j} \rangle = {\delta_{ij}} i !, \quad \rho(z) = \frac{1}{\sqrt{2\pi}}\exp^{-z^{2}/{2}}.
\end{align*}
Subsequently, the parameters, e.g. $d\sim \pazocal{N}(\mu_d,\sigma^{2}_{d})$, can be expanded by the standard random variable $\xi_{d} \sim \pazocal{N}(0,1)$, which will be an element of the random vector $\xi\in \mathbb{R}^{3}$, i.e.:
\begin{align*}
    d(\xi_{d}) = \mu_{d}\phi_{0}(\xi_{d}) + \sigma_{d}\phi_{1}(\xi_{d}) = \mu_{d} + \sigma_{d}\xi_{d}.
\end{align*}
%
The dynamics of the robot are provided by: 
\begin{align*}
\begin{bmatrix}
x_{k+1} \\
y_{k+1} \\
\theta_{k+1} 
\end{bmatrix}
=
\begin{bmatrix}
x_{k}\\
 y_{k}\\
 \theta_{k}
\end{bmatrix}
+
\begin{bmatrix}
v \cos{\theta_{k}}\\
 v \sin{\theta_{k}}\\
 \omega
\end{bmatrix}dt,
\end{align*}
where $x,y$ are its position coordinates, $\theta$ is its orientation, $v$ and $\omega$ are its translational and rotational velocities, respectively, and $dt$ is the time step.
The control inputs of the robot are the rotational velocities of its two wheels: $u_1 = \dot{\phi}_{R}$ and $u_2 = \dot{\phi}_{L}$. The translational velocities of the wheels can be obtained from:
\begin{align*}
v_{R} = 
r_{R}u_{1}, \quad
v_{L} = 
r_{L}u_{2},
\end{align*}
while $v$ and $\omega$ are given by: 
\begin{align*}
 v = \frac{1}{2}(v_{R}+v_{L}), \quad
 \omega = \frac{v_{R}-v_{L}}{2d}.
\end{align*}
In this setting, chance constraints can be formulated by imposing the probability that the robot will avoid colliding with an obstacle to be greater than $p$.
The confidence region of the robot can be used to transform the chance constraints into deterministic ones. In particular, we approximate the area $R$ where the robot might exist with probability $p$ by using the ellipsoid:
\[ R = \Set{z = (x,y)^{\mathsf{T}} |
\begin{array}{l}
 (z-\bar{z})^{\mathsf{T}}\Sigma^{-1}(z-\bar{z}) \leq s(p)
\end{array}}, \]
where $\bar{z}$ is mean of $z$, $\Sigma$ is covariance matrix given by:
\begin{align*}
\Sigma =
    \begin{bmatrix}
    \mathbb{VAR}[x] & \mathbb{COV}[xy]\\
    \mathbb{COV}[xy] & \mathbb{VAR}[y]
    \end{bmatrix},
\end{align*}
and $s(p)$ is the scaling factor. If the position of the robot followed a Gaussian distribution, then $s(p)$ could be obtained by computing the chi-square inverse cumulative distribution function with a specified $p$, which is fixed over time. In our case, however, the distribution of the position is not Gaussian even though the uncertain parameters follow such a distribution. Thus, we evaluate $s(p)$ using Monte Carlo simulation by sampling the random parameters at each time step. In order to facilitate the computation, we overestimate the confidence region with a circle whose radius is the major axis of the ellipsoid, which can be obtained as the largest eigenvalue of $\Sigma$, i.e. $\lambda_{\rm{max}}$.
Thus, a constraint for a circle obstacle with center $(x_{c}, y_{c})$ and radius $r_{c}$ can be formulated as follows:
\begin{align*}
g(x,y) =& (r_{c} + \sqrt{s(p)\lambda_{\rm{max}}})^2 \\
-& [(x-x_{c})^{2} + (y-y_{c})^{2}] \leq 0
\end{align*}
Note that the gradient and the Hessian of $g$ - which are required for constrained DDP - can be analytically computed with the derivative of $\lambda_{\rm{max}}$ \cite{Nelson1976eigenvectorderivative}.
\begin{figure}
\centering
\hspace*{-1em}
    \subfigure[\label{fig:robot_traj_gpc-cddp}]
    {
       \includegraphics[trim={0.1cm 0.1cm 0.4cm 0.4cm},clip,width=0.46\columnwidth]{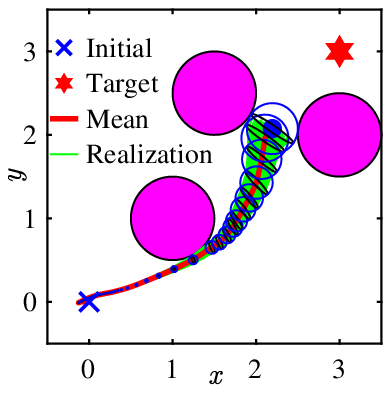}
    }\hspace*{-0.5em}
    \subfigure[\label{fig:robot_traj_gpc-mpc-cddp}]
    {
        \includegraphics[trim={0.1cm 0.1cm 0.4cm 0.4cm},clip,width=0.46\columnwidth]{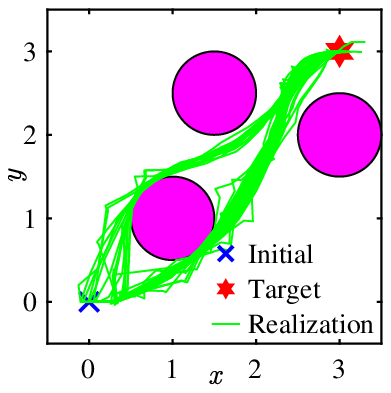}
    }\hspace*{-0.5em}\\
    \hspace*{-1em}
    \subfigure[\label{fig:robot_traj_gpc-mpc-cddp_mpc_call}]
    {
        \includegraphics[trim={0.1cm 0.1cm 0.4cm 0.4cm},clip,width=0.46\columnwidth]{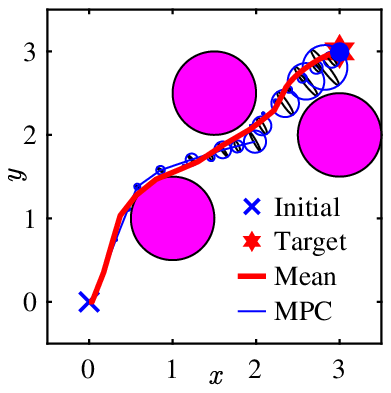}
    }\hspace*{-0.5em}
        \subfigure[\label{fig:robot_normal_mpc}]
    {
        \includegraphics[trim={0.1cm 0.1cm 0.4cm 0.4cm},clip,width=0.46\columnwidth]{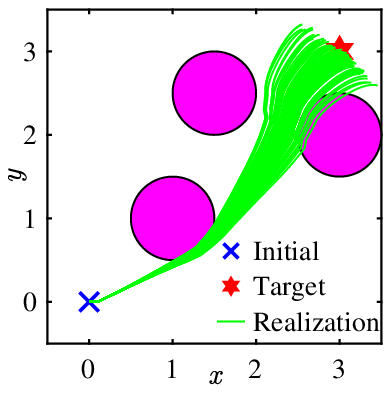}
    }
    \caption{Differential wheeled robot: Reaching task from  $[0,0,0]^{\mathsf{T}}$,  to  $[3,3,0]^{\mathsf{T}}$ while  avoiding three obstacles (pink). Obtained trajectories from
    (a) gPC CDDP: mean variance, and 100 realizations, (b) MPC gPC DDP: 100 realizations, (c) MPC gPC DDP: one realization, mean and variance from several MPC cycles, (d) MPC DDP: 100 realizations.}
  \label{fig:robot_traj} 
\end{figure}
In this simulation, we have set $dt = 0.02$s, $N = 60$, $H = 10$, $r=2$, and $p =0.95$. The initial control is set to be zero and the control bounds to be $|u_i|\leq 100$, $i=1,2$.

In Fig. \ref{fig:robot_traj}, we demonstrate a comparison between gPC CDDP, gPC MPC CDDP and MPC CDDP which is purely deterministic. In Fig. \ref{fig:robot_traj_gpc-cddp}, the accumulation of the variance while using gPC CDDP is demonstrated, making the robot unable to pass through the two obstacles. The confidence region is validated by overlaying the realizations of the trajectories with sampled parameters. In Fig. \ref{fig:robot_traj_gpc-mpc-cddp},
gPC MPC CDDP is applied to the dynamics obtained by the realizations of different sampled parameters. The results show that the method can successfully reduce the variance growth, so the robot reaches the desired target. Note that for this method, the exact values of the dynamics parameters are not required, but only information of the states at every time step. Fig. \ref{fig:robot_traj_gpc-mpc-cddp_mpc_call} shows one trajectory from the same algorithm and the confidence region from some MPC cycles. Note that the trajectory indeed lies inside the confidence region computed by the algorithm. Finally, the performance of MPC CDDP with the parameters being sampled, is presented in Fig \ref{fig:robot_normal_mpc}. In this simulation, the parameters used in the model are the mean parameters. Unlike gPC CDDP, the uncertainty of the dynamics which arises from the parameters is not considered here. Therefore, the real dynamics trajectory may violate the constraints even though they appear to be satisfied by the model trajectories computed by MPC. The poor performance of this method validates that not only using MPC but also taking into consideration the uncertainty is necessary.

\subsection{Quadrotor}
\begin{figure}
\centering
\hspace*{-1.3em}
    \subfigure[\label{fig:quad_traj_gpc-cddp}]
    {
       \includegraphics[clip, trim={1.02cm 0.1cm 0.6cm 0.1cm},width=0.55\columnwidth]{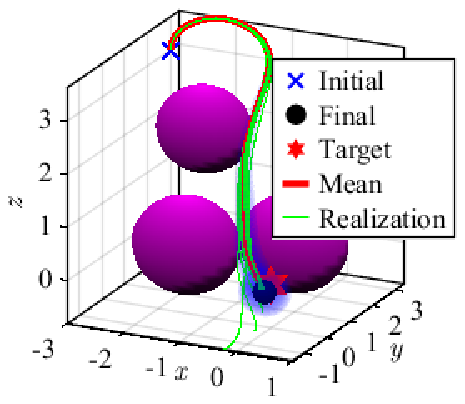}

    }
    \hspace*{-1.6em}
    \subfigure[\label{fig:quad_traj_gpc-mpc-cddp}]
    {
    \includegraphics[clip,trim={0.1cm 0.1cm 0.25cm 0.2cm},width=0.47\columnwidth]{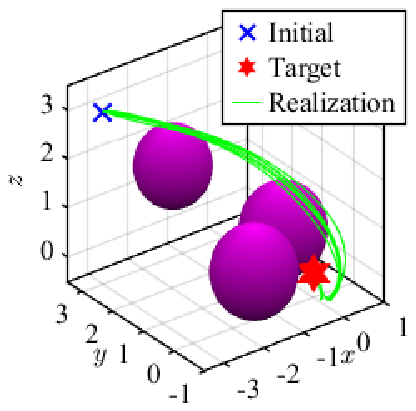}
    }
    \caption{Quadrotor: Reaching task from  $[3,-3,3]^{\mathsf{T}}$  to  $[0,0,0]^{\mathsf{T}}$ and staying at the target while dodging three obstacles (pink). Obtained trajectories from
    (a) gPC CDDP: mean variance (blue ball), and $100$ realizations, (b) gPC MPC CDDP: 25 realizations.}
  \label{fig:quad_traj} 
\end{figure}

Subsequently, we test our method on a quadrotor whose dynamics are provided in \cite{Luukkonen2011}. We assume that we are aware of the upper and lower limits of the drag and lift coefficients. Thus, they can be viewed as uniformly distributed random variables which can be described using Legendre polynomials - according to Table I. 
The first few polynomials, their inner products, and the PDF are as follows:
\begin{align*}
   &\phi_{0}(z) = 1,\quad \phi_{1}(z) = z, \quad \phi_{2}(z) = \frac{3}{2}z^{2}-\frac{1}{2},\dotsc, \\ 
    &\langle \phi_{i}, \phi_{j} \rangle = \frac{1}{2i+1}{\delta_{ij}}, \quad
    \rho(z) = \frac{1}{2}.
\end{align*}
with $z \in \mathbb{R}$ and $i,j = 0,1, \dotsc$.
Note that the inner product is scaled by $1/2$ compared to the orthogonality relation of the polynomial\cite{Xiu2010gpc_text}. Each uncertain parameter, e.g. the drag coefficient $k_{d}$, is expanded by its mean $\mu_{k_{d}}$ and the distance between its mean and lower or upper limit of the distribution denoted by $\Delta_{k_{d}}$, i.e.:
\begin{align*}
k_{d}(\xi_{1}) = \mu_{k_{d}}\phi_{0}(\xi_{1})+ \Delta_{k_{d}}\phi_{1}(\xi_{1}) = \mu_{k_{d}}+ \Delta_{k_{d}}\xi_{1},  
\end{align*}
where $\xi_{1}$ is a uniformly distributed scalar random variable in $[-1, 1]$. Similarly, the lift coefficient $k_{l}$ can be expanded with the standard uniformly distributed random variable $\xi_{2}$. Thus, we have the random vector $\xi = [\xi_{1}, \xi_{2}]^{\mathsf{T}}$.
We have used $\mu_{k_{b}} = 1.140\times10^{-7}$, $\Delta_{b} = \mu_{k_{b}}/3$,
$\mu_{k_{l}} = 2.980\times10^{-6}$, $\Delta_{k_{l}} = \mu_{k_{l}}/3$.
In this simulation, the control of the system is the force generated by four rotors with control limits $0 \leq u\leq 3$. The initial state of the quadrotor is hovering. Moreover, $N = 100$, $H = 25$, $dt=0.02$s, $r = 2$, and $p=0.95$. Although this problem is 3D, the concepts of confidence region and chance constraints presented in \ref{sec:result_robot} are readily extended and applied to this task.

Fig. \ref{fig:quad_traj} demonstrates the quadrotor simulation results. In Fig. 
\ref{fig:quad_traj_gpc-cddp}, we observe that most of the realizations of the trajectories are captured within the blue ball which shows the confidence region of the quadrotor. We also observe again that with gPC CDDP, the growth of the variance pushes the quadrotor away from the desired state in order to satisfy the constraints. 
The trajectories obtained from gPC MPC CDDP are presented in Fig. \ref{fig:quad_traj_gpc-mpc-cddp}. Similarly to the case of the differential wheeled robot, our receding horizon approach successfully reduces the variance, driving the trajectories closer to the target.  

\section{Conclusion} 
\label{sec:conclusion}
In this work, we propose a novel chance-constrained receding horizon control method that is able to handle uncertainty arising from model parameters using gPC.
In particular, by transforming a chance-constrained optimization problem into a deterministic one, we are able to solve it using a constrained DDP technique. 
Our method successfully handles high-dimensional dynamics while enjoying the convergence properties and scalability of DDP. Furthermore, in order to address the accumulation of uncertainty on the trajectory, we implement our algorithm in a receding horizon fashion. Simulation results demonstrate that our method is applicable to complex robotic tasks and able to deal effectively with parametric uncertainty and nonlinear constraints.


Future work would examine the robustness of our method under other forms of uncertainty such as process noise. 
Moreover, we would like to explore potential combinations of our method with others that can handle non-parametric uncertainty.
Finally, we are interested in incorporating learning techniques such as Gaussian processes that will lead to safe learning algorithms under parametric and non-parametric uncertainty.


\section*{Acknowledgments}
The work of the second and third authors was supported by awards NSF CMMI-1936079 and NSF CPS-1932288, respectively. Augustinos Saravanos also acknowledges support by the A. Onassis Foundation Scholarship.


\end{document}